\documentclass[11pt]{amsart}
\input{diagrams}

\usepackage{amsmath}
\usepackage{amsthm}
\usepackage{amssymb}
\usepackage{amsfonts}
\usepackage{amsxtra}
\usepackage{amscd}
\usepackage{epsfig}
\usepackage{verbatim}
\usepackage{latexsym,amstext,epsfig}

\usepackage[all, knot]{xy}
\xyoption{arc}

\newsymbol\pp 1275
\newcommand{\Hom}{\operatorname{Hom}}
\newcommand{\End}{\operatorname{End}}

\newcommand{\Rep}{\operatorname{Rep}}
\newcommand{\Proj}{\operatorname{Proj}}
\newcommand{\SI}{\operatorname{SI}}
\newcommand{\SL}{\operatorname{SL}}
\newcommand{\GL}{\operatorname{GL}}
\newcommand{\PGL}{\operatorname{PGL}}
\newcommand{\ZZ}{\mathbb Z}
\newcommand{\CC}{\mathbb C}

\newcommand{\NN}{\mathbb N}
\newcommand{\QQ}{\mathbb Q}

\newcommand{\Id}{\operatorname{Id}}

\newtheorem{theorem}{Theorem}[section]
\newtheorem{proposition}[theorem]{Proposition}

\newtheorem{lemma}[theorem]{Lemma}

\theoremstyle{definition}
\newtheorem{definition}[theorem]{Definition}
\newtheorem{remark}[theorem]{Remark}
\newtheorem{conjecture}[theorem]{Conjecture}
\newtheorem{example}[theorem]{Example}

\newcount\cols
{\catcode`,=\active\catcode`|=\active
\gdef\Young(#1){\hbox{$\vcenter
{\mathcode`,="8000\mathcode`|="8000
\def,{\global\advance\cols by 1 &}%
\def|{\cr
      \multispan{\the\cols}\hrulefill\cr
       &\global\cols=2 }%
  \offinterlineskip\everycr{}\tabskip=0pt
  \dimen0=\ht\strutbox \advance\dimen0 by \dp\strutbox
    \halign
    {\vrule height \ht\strutbox depth \dp\strutbox##
      &&\hbox to \dimen0{\hss$##$\hss}\vrule\cr
     \noalign{\hrule}&\global\cols=2 #1\crcr
     \multispan{\the\cols}\hrulefill\cr%
   }
}$}} }

\title[]{Counterexamples to Okounkov's log-concavity conjecture}

\author{Calin Chindris}
\address{University of Minnesota, School of Mathematics, Minneapolis, MN, USA}
\email[Calin Chindris]{chindris@math.umn.edu}

\author{Harm Derksen}
\address{University of Michigan, Department of Mathematics, Ann Arbor, MI, USA}
\email[Harm Derksen]{hderksen@umich.edu}

\author{Jerzy Weyman}
\address{Northeastern University, Department of Mathematics, Boston, MA, USA}
\email[Jerzy Weyman]{j.weyman@neu.edu}
\thanks{The second author was supported by NSF, grant DMS 0349019 and the third author was supported by NSF, grant 0600229}

\date{\today}

\begin{document}
\bibliographystyle{plain}
 %\keywords{...}
\begin{abstract}
We give counterexamples to Okounkov's log-concavity conjecture for
Littlewood-Richardson coefficients.
\end{abstract}
\maketitle

\section{Introduction}

Motivated by physical considerations, Okounkov \cite[Conjecture
1]{Oko} has conjectured that the Littlewood-Richardson
coefficients $c_{\mu, \nu}^{\lambda}$ are log-concave in
$(\lambda, \mu, \nu).$ A particular version of this conjecture
would be (see also \cite[pp. 239]{F1}):

\begin{conjecture}[$\textbf{Okounkov's log-concavity conjecture}$]\label{conj-log} Let $\lambda, \mu, \nu$ be three
partitions. Then
$$
c_{(N+1)\mu, (N+1)\nu}^{(N+1)\lambda} \cdot c_{(N-1)\mu,
(N-1)\nu}^{(N-1)\lambda} \leq (c_{N\mu, N\nu}^{N\lambda})^2,
$$
for every integer $N \geq 1.$
\end{conjecture}

Important implications of this conjecture are also discussed in
\cite{Oko}. It is easy to see that Conjecture \ref{conj-log}, if
true, would immediately imply a conjecture of Fulton on
Littlewood-Richardson coefficients (see \cite{Bel} or \cite{KTW}).
Moreover, the log-concavity of the Littlewood-Richardson
coefficients as a function of highest weights would imply the
saturation conjecture for Littlewood-Richardson coefficients (see
\cite{DW1} or \cite{KT}) and the Schur-log concavity conjecture
for skew-Schur functions (see \cite{LPP}). We should also point
out that the results in \cite{TZ} give some evidence for the
log-concavity conjecture. However, the conjecture turns out to be
false in general.

\smallskip
In this note, we construct infinite families of counterexamples to
Conjecture \ref{conj-log} (and hence to the original Okounkov's
conjecture):

\begin{theorem}\label{prop-cex} Let $n \geq 1$ be an integer and let
$\lambda(n),$ $\mu(n)$  be two partitions defined by
$$
\lambda(n)=(4^n,3^{2n},2^n),\text{~and~} \mu(n)=(3^n,2^n,1^n).
$$
Then

$$
c_{\mu(n),\mu(n)}^{\lambda(n)}=\binom{n+2}{2} \text{~and~}
c_{2\mu(n),2\mu(n)}^{2\lambda(n)}= \binom{n+5}{5}.
$$
Consequently, when $n\geq 21$, Conjecture \ref{conj-log} fails for
$\lambda=\lambda(n),$ $\mu=\nu=\mu(n),$ and $N=1.$
\end{theorem}

The details of our notations can be found in the notation
paragraph at the end of this section. We would like to point out
that Conjecture \ref{conj-log} is true \emph{asymptotically}. This
fact was proved by Okounkov in \cite[Section 3.5]{Oko}. A
different proof can be found in Section \ref{poly-sec} (Remark
\ref{asympt-rmk} and Example \ref{asympt-ex}).

The layout of this note is as follows. In Section
\ref{LRproof-sec}, we give a direct proof of Theorem
\ref{prop-cex} by using the Littlewood-Richardson rule. A
different approach to Okounkov's conjecture is based on quiver
theory. In Section \ref{quiver-sec}, we review some tools from
quiver invariant theory and explain why the log-concavity
conjecture is bound to fail (see Section \ref{poly-sec}). In
Section \ref{cex-sec}, we give another proof of Theorem
\ref{prop-cex} and present more counterexamples. In particular,
Proposition \ref{parK-cex-prop} provides counterexamples to the
log-concavity conjecture \cite[Conjecture 6.17]{Kir} for
\emph{parabolic} Kostka numbers.

\subsection*{Notations} A partition is a sequence $\lambda=(\lambda_1, \dots,
\lambda_r)$ of integers such that $\lambda_1 \geq \dots \lambda_r
\geq 0$. The length of a partition is defined to be the number of
its non-zero parts. If $\lambda$ is a partition, we define
$|\lambda|$ to be the sum of its parts. The Young diagram of a
partition $\lambda=(\lambda_1, \dots, \lambda_r)$ is a collection
of boxes, arranged in left-justified rows with $\lambda_i$ boxes
in row $i$. For a partition $\lambda,$ we denote by $\lambda'$ the
partition conjugate to $\lambda$, i.e., the Young diagram of
$\lambda'$ is the Young diagram of $\lambda$ reflected with
respect to its main diagonal.

If $\lambda=(\lambda_1, \dots, \lambda_r)$ is a partition then we
define $N\lambda$ by $N\lambda=(N\lambda_1, \dots, N\lambda_r).$
By $\lambda=(\lambda_1^{m_1},\dots,\lambda_k^{m_k}),$ we denote
the partition that has $m_i$ parts equal to $\lambda_i,$ $1 \leq i
\leq k.$ For a partition $\lambda$ of length at most $r$,
$S^{\lambda}(V)$ denotes the irreducible polynomial representation
of $\GL(V)$ with highest weight $\lambda,$ where $V$ is an
$r$-dimensional complex vector space. Let $\lambda, \mu, \nu$ be
three partitions of length at most $r$. Then we define the
Littlewood-Richardson coefficient $c_{\mu,\nu}^{\lambda}$ to be
the multiplicity of $S^{\lambda}(V)$ in $S^{\mu}(V) \otimes
S^{\nu}(V)$, i.e.,
$$
c_{\mu,\nu}^{\lambda}=\dim_{\CC}(S^{\lambda}(V)^* \otimes
S^{\mu}(V) \otimes S^{\nu}(V))^{\GL(V)},
$$
where $S^{\lambda}(V)^*$ is the dual representation.  More
generally, if $\gamma, \lambda(1), \dots, \lambda(m)$ are
partitions of length at most $r,$ we define
$$
c^{\gamma}_{\lambda(1), \dots, \lambda(m)}=
\dim_{\mathbb C} (S^{\gamma}(V)^{*}\otimes S^{\lambda(1)}(V)
\otimes \dots \otimes S^{\lambda(m)}(V))^{\GL(V)}.
$$

\section{A direct proof by Littlewood-Richardson rule} \label{LRproof-sec}

Our main references for Young tableau and Littlewood-Richardson
rule are \cite{F3} and \cite{Mac} (see also \cite{F1}). If
$\lambda ,\mu ,\nu$ are three partitions, the
Littlewood-Richardson coefficient $c^{\lambda}_{\mu, \nu}$ can be
described as the cardinality of the set $LR(\lambda ,\mu ,\nu)$ of
diagrams $D$ of skew shape $\lambda /\mu$, filled with $\nu_1$
$1$'s, $\nu_2$ $2$'s, etc., subject to the following conditions:
\begin{enumerate}
\renewcommand{\theenumi}{\arabic{enumi}}
\item $D$ is a semistandard Young tableau, i.e., the entries in
rows are weakly increasing from left to right and the entries in
columns are strictly increasing from top to bottom;

\item $D$ is a lattice permutation, i.e., when the entries are
listed, from right to left in rows, starting with the top row, the
resulting word $w(D)$ is a lattice permutation. This last
condition means that for any integer $1 \leq r \leq |\nu|,$ and
any positive integer $i$, the number of occurrences of $i$ in the
first $r$ entries of $w(D)$ is no less than the number of
occurrences of $i+1$ in these first $r$ entries.
\end{enumerate}

\begin{example} For $\lambda=(4,2,1),$  $\mu=(3,1,0),$ and $\nu=
(2,1,0),$ there are only two diagrams $D$ that satisfy conditions
$(1)$ and $(2)$ above:
$$\Young (,,,1|,1|2), \Young (,,,1|,2|1).$$
Note that the diagram
$$\Young (,,,2|,1|1)$$
is not a lattice permutation. Reversing the roles of $\mu$, $\nu,$
we get the diagrams
$$\Young (,,1,1|,1|2), \Young (,,1,1|,2|1).$$
\end{example}

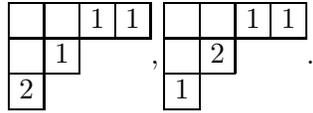
\begin{proof}[Proof of Theorem \ref{prop-cex}]
For our purposes it will be convenient to work with conjugate
partitions. This is always possible since the
Littlewood-Richardson coefficients are invariant when passing to
conjugate partitions. First, we show that
$$c^{\lambda(n)'}_{\mu(n)',\mu(n)'} = {{n+2}\choose 2},$$
where $\lambda(n)'=(4n,4n,3n,n)$ and $\mu(n)'=(3n,2n,n)$. We look
at the cases $n=1,2$ and then describe the general pattern. For
$n=1$, the multiplicity is 3 because the only three skew diagrams
satisfying the requirements of the Littlewood-Richardson rule are:
$$
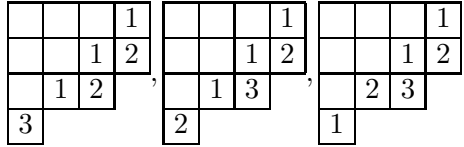
\Young (,,,1|,,1,2|,1,2|3) ,\Young (,,,1|,,1,2|,1,3|2) ,\Young
(,,,1|,,1,2|,2,3|1)
$$

Let us call these tableaux $S_3$, $S_2$, $S_1$ respectively, i.e.,
we label them by the content of the last row.

For $n=2$, we get six tableaux
$$\Young (,,,,,,1,1|,,,,1,1,2,2|,,1,1,2,2|3,3), \Young (,,,,,,1,1|,,,,1,1,2,2|,,1,1,2,3|2,3) ,\Young (,,,,,,1,1|,,,,1,1,2,2|,,1,2,2,3|1,3)
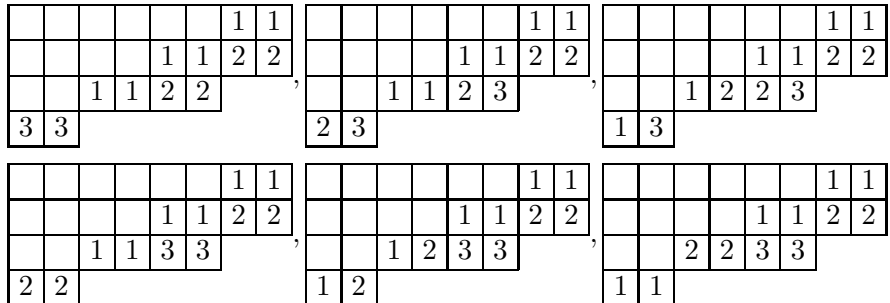
$$
$$\Young (,,,,,,1,1|,,,,1,1,2,2|,,1,1,3,3|2,2), \Young (,,,,,,1,1|,,,,1,1,2,2|,,1,2,3,3|1,2) ,\Young (,,,,,,1,1|,,,,1,1,2,2|,,2,2,3,3|1,1)
$$ which, looking at their last row, clearly correspond to the
monomials of degree $2$ in $S_1, S_2 ,S_3$.

So, for general $n$ we define a bijection between the set of
monomials
$$S^a := S_1^{a_1} S_2^{a_2} S_3^{a_3}$$
of degree $n$ and the set $LR(\lambda(n)', \mu(n)', \mu(n)')$ of
tableaux of the shape $\lambda(n)'/\mu(n)'$ satisfying the
conditions $(1),$ $(2)$ above, whose cardinality is equal to
$c^{\lambda(n)'}_{\mu(n)',\mu(n)'}$. To achieve this, we associate
to each monomial $S^a$ a tableau $E(a)$ from $LR(\lambda(n)',
\mu(n)', \mu(n)')$ as follows. The first two rows of each $E(a)$
are the same, they contain in each column the numbers $1$, for
columns with numbers $2n+1,\ldots ,3n$, and $1,2$ for columns with
numbers $3n+1 ,\ldots ,4n$.

The filling of the last row of $E(a)$ is a tableau of shape $(n)$
and we just define it to have $a_1$ $1$'s, $a_2$ $2$'s and $a_3$
$3$'s. Now it is clear that the remaining third row of
$\lambda(n)'/\mu(n)'$ can be uniquely filled by the remaining
available numbers to get the tableau $E(a)$ from $LR(\lambda(n)',
\mu(n)', \mu(n)')$. Indeed, the first (from left to right) $n$
boxes in the third row have to be filled by the remaining $1$'s
and $2$'s and the last $n$ boxes by the remaining $2$'s and $3$'s,
in weakly increasing order. This assures semi-standardness. The
lattice permutation condition between $1$'s and $2$'s is satisfied
because there are $2n$ $1$'s already in the columns $2n+1 ,\ldots
,4n.$ The lattice permutation condition between $2$'s and $3$'s is
also satisfied because there are already $n$ $2$'s in columns
$3n+1,\ldots ,4n.$

This gives us an injection from the set of monomials $S^a$ to the
set of tableaux $LR(\lambda(n)', \mu(n)', \mu(n)')$. It is clearly
surjective because to each diagram $E$ from  the set
$LR(\lambda(n)', \mu(n)', \mu(n)')$ we can associate the monomial
$S_1^{u_1}S_2^{u_2} S_3^{u_3}$ of degree $n$ by taking $u_i$ to be
the number of occurrences of $i$ in the last row of $E$. This
shows that $$c^{\lambda(n)'}_{\mu(n)',\mu(n)'} ={{n+2}\choose
2}.$$

Let us turn to the second statement. We need to show that
$$c^{\rho(n)}_{\sigma(n),\sigma(n)} = {{n+5}\choose 5},$$
where $\rho(n)=(2\lambda(n))' =(4n,4n,4n,4n,3n,3n,n,n)$ and
$\sigma(n)=(2\mu(n))'=(3n,3n,2n,2n,n,n)$. Let us exhibit the case
$n=1:$
$$
\Young (,,,1|,,,2|,,1,3|,,2,4|,1,3|,2,4|5|6) , \Young
(,,,1|,,,2|,,1,3|,,2,4|,1,3|,2,5|4|6), \Young
(,,,1|,,,2|,,1,3|,,2,4|,1,5|,2,6|3|4)
$$
$$
\Young (,,,1|,,,2|,,1,3|,,2,4|,1,3|,4,5|2|6) , \Young
(,,,1|,,,2|,,1,3|,,2,4|,1,5|,3,6|2|4),\Young
(,,,1|,,,2|,,1,3|,,2,4|,3,5|,4,6|1|2)
$$

Let us label these tableaux by the content of the first column,
i.e., $T_{5,6}$,  $T_{4,6}$, $T_{3,4}$, $T_{2,6}$, $T_{2,4}$, and
$T_{1,2}$, respectively. We also order them by a total order
respecting the lexicographic order of the indices, i.e.,
$$T_{1,2}<T_{2,4}<T_{2,6}<T_{3,4}<T_{4,6}<T_{5,6} .$$
We define a bijection between the set of monomials
$$T^a := T_{1,2}^{a_{1,2}} T_{2,4}^{a_{2,4} }T_{2,6}^{a_{2,6}} T_{3,4}^{a_{3,4}} T_{4,6}^{a_{4,6}} T_{5,6}^{a_{5,6}}$$
of degree $n$ and the set $LR(\rho(n), \sigma(n), \sigma(n))$ of
tableaux of the shape $\rho(n)/\sigma(n)$ satisfying the
conditions $(1),$ $(2)$ above, whose cardinality is equal to
$c^{\rho(n)}_{\sigma(n),\sigma(n)}$. To achieve this, we associate
to each monomial $T^a$ a tableau $D(a)$ from $LR(\rho(n),
\sigma(n), \sigma(n))$ as follows. The first four rows of each
$D(a)$ are the same, they contain in each column the numbers
$1,2$, for columns with numbers $2n+1,\ldots ,3n$, and $1,2,3,4$
for columns with numbers $3n+1 ,\ldots ,4n.$

The filling of the last two rows of $D(a)$ form a tableau of shape
$(n^2)$. We start with a tableau having $a_{i,j}$ columns of type
$\Young(i|j)$. We order them according to the order on $T_{i,j}$,
so columns $\Young (1|2)$ are the first ones and the columns
$\Young (5|6)$ are the last ones. The only problem is that the
columns $\Young(3|4)$ and $\Young (2|6)$ cannot be standard in any
order. So, every occurrence of the columns $\Young (3|4)$ and
$\Young (2|6)$ has to be replaced by $\Young (2,3|4,6)$. This
defines the filling the last two rows of $D(a)$. Now, we claim
that the remaining fifth and sixth rows of $\rho(n)/\sigma(n)$ can
be uniquely filled by the remaining available numbers to get the
tableau $D(a)$ from $LR(\rho(n), \sigma(n), \sigma(n)).$ The point
is that $1$'s have to appear in the fifth row at the beginning,
and $2$'s cannot appear in that row after $1$'s because lattice
permutation condition would be violated. Similarly, $6$'s have to
appear at the end of the sixth row, but there has to be a $5$
above each $6$, otherwise lattice permutation condition is
violated. The rest of $5$'s have to appear before the $6$'s. The
remaining part of the diagram can be uniquely filled with $3$'s
and $4$'s to complete it to a standard diagram. Indeed, the number
$4$ cannot appear in the fifth row, because the lattice
permutation condition would be violated (the number of $3$'s and
$4$'s in the first four rows is the same). Semi-standardness and
the lattice permutation condition easily follows. This gives us an
injection from the set of monomials $T^a$ to the set $LR(\rho(n),
\sigma(n), \sigma(n))$. This is enough for the counterexample,
because we showed that the coefficient $c^{\rho(n)}_{ \sigma(n),
\sigma(n)} $ is at least ${n+5}\choose 5$. The fact that the
defined map is surjective is not difficult to prove, so we leave
it to the reader.
\end{proof}

\section{Quiver theory} \label{quiver-sec} In this section we
review the main tools from quiver invariant theory that will be
used to study Littlewood-Richardson coefficients.
\subsection{Generalities} A quiver $Q=(Q_0,Q_1,t,h)$
consists of a finite set of vertices $Q_0$, a finite set of arrows
$Q_1$ and two functions $t,h:Q_1 \to Q_0$ that assign to each
arrow $a$ its tail $ta$ and its head $ha,$ respectively. We write
$ta{\buildrel a\over\longrightarrow}ha$ for each arrow $a \in
Q_1$.

For simplicity, we will be working over the field of complex
numbers $\CC.$ A representation $V$ of $Q$ over $\CC$ is a family
of finite dimensional $\CC$-vector spaces $\lbrace V(x) \mid x\in
Q_0\rbrace$ together with a family $\{ V(a):V(ta)\rightarrow V(ha)
\mid a \in Q_1 \}$ of $\CC$-linear maps. If $V$ is a
representation of $Q$, we define its dimension vector $\underline
d_V$ by $\underline d_V(x)=\dim_{\CC} V(x)$ for every $x\in Q_0$.
Thus the dimension vectors of representations of $Q$ lie in
$\Gamma=\ZZ^{Q_0}$, the set of all integer-valued functions on
$Q_0$. For every vertex $x$, we denote by $e_x$ the simple
dimension vector corresponding to $x$, i.e. $e_x(y)=\delta_{x,y},
\forall y \in Q_0,$ where $\delta_{x,y}$ is the Kronecker symbol.

Given two representations $V$ and $W$ of $Q$, we define a morphism
$\phi:V \rightarrow W$ to be a collection of linear maps $\lbrace
\phi(x):V(x)\rightarrow W(x)\mid x \in Q_0 \rbrace$ such that for
every arrow $a\in Q_1$, we have $\phi(ha)V(a)=W(a)\phi(ta)$. We
denote by $\Hom_Q(V,W)$ the $\CC$-vector space of all morphisms
from $V$ to $W$. In this way, we obtain the abelian category
$\Rep(Q)$ of all quiver representations of $Q.$ Let $V$ and $W$ be
two representations of $Q.$ We say that $V$ is a subrepresentation
of $W$ if $V(x)$ is a subspace of $W(x)$ for all vertices $x \in
Q_0$ and $V(a)$ is the restriction of $W(a)$ to $V(ta)$ for all
arrows $a \in Q_1.$

If $\alpha,\beta$ are two elements of $\Gamma$, we define the
Euler inner product
\begin{equation}
\langle\alpha,\beta \rangle = \sum_{x \in Q_0}
\alpha(x)\beta(x)-\sum_{a \in Q_1} \alpha(ta)\beta(ha).
\end{equation}

\emph{From now on, we will assume that our quivers are without
oriented cycles.}

\subsection{Semi-invariants for quivers}
Let $\beta$ be a dimension vector of $Q$. The representation space
of $\beta-$dimensional representations of $Q$ is defined by
$$\Rep(Q,\beta)=\bigoplus_{a\in Q_1}\Hom(\CC^{\beta(ta)}, \CC^{\beta(ha)}).$$
If $\GL(\beta)=\prod_{x\in Q_0}\GL(\beta(x))$ then $\GL(\beta)$
acts algebraically on $\Rep(Q,\beta)$ by simultaneous conjugation,
i.e., for $g=(g(x))_{x\in Q_0}\in \GL(\beta)$ and $V=\{V(a)\}_{a
\in Q_1} \in \Rep(Q,\beta),$ we define $g \cdot V$ by
$$(g\cdot V)(a)=g(ha)V(a)g(ta)^{-1}\ \text{for each}\ a \in Q_1.$$ In
this way, $\Rep(Q,\beta)$ is a rational representation of the
linearly reductive group $\GL(\beta)$ and the $\GL(\beta)-$orbits
in $\Rep(Q,\beta)$ are in one-to-one correspondence with the
isomorphism classes of $\beta-$dimensional representations of $Q.$
As $Q$ is a quiver without oriented cycles, one can show that
there is only one closed $\GL(\beta)-$orbit in $\Rep(Q,\beta)$ and
hence the invariant ring $\text{I}(Q,\beta)= \CC
[\Rep(Q,\beta)]^{\GL(\beta)}$ is exactly the base field $\CC.$

Now, consider the subgroup $\SL(\beta) \subseteq \GL(\beta)$
defined by
$$
\SL(\beta)=\prod_{x \in Q_0}\SL(\beta(x)).
$$

Although there are only constant $\GL(\beta)-$invariant polynomial
functions on $\Rep(Q,\beta)$, the action of $\SL(\beta)$ on
$\Rep(Q,\beta)$ provides us with a highly non-trivial ring of
semi-invariants. Note that any $\sigma \in \ZZ^{Q_0}$ defines a
rational character of $\GL(\beta)$ by
$$\{g(x) \mid x \in Q_0\} \in \GL(\beta) \mapsto \prod_{x \in
Q_0}(\det g(x))^{\sigma(x)}.$$ In this way, we can identify
$\Gamma=\ZZ ^{Q_0}$ with the group $X^\star(\GL(\beta))$ of
rational characters of $\GL(\beta),$ assuming that $\beta$ is a
sincere dimension vector (i.e. $\beta(x)>0$ for all vertices $x
\in Q_0$). We also refer to the rational characters of
$\GL(\beta)$ as weights.

Let $\SI(Q,\beta)= \CC [\Rep(Q,\beta)]^{\SL(\beta)}$ be the ring
of semi-invariants. As $\SL(\beta)$ is the commutator subgroup of
$\GL(\beta)$ and $\GL(\beta)$ is linearly reductive, we have
$$\SI(Q,\beta)=\bigoplus_{\sigma
\in X^\star(\GL(\beta))}\SI(Q,\beta)_{\sigma},
$$
where $$\SI(Q,\beta)_{\sigma}=\lbrace f \in \CC [\Rep(Q,\beta)]
\mid gf= \sigma(g)f \text{~for all~}g \in \GL(\beta)\rbrace$$ is
the space of semi-invariants of weight $\sigma.$ If $\alpha \in
\Gamma$, we define $\sigma=\langle \alpha,\cdot \rangle$ by
$$\sigma(x)=\langle \alpha,e_x \rangle, \forall x\in
Q_0.$$ Similarly, one can define the weight $\tau=\langle
\cdot,\alpha \rangle.$

\begin{lemma}[Reciprocity Property]\cite[Corollary 1]{DW1}\label{reciprocity}
Let $\alpha$ and $\beta$ be two dimension vectors. Then
$$\dim\SI(Q,\beta)_{\langle \alpha, \cdot \rangle} = \dim
\SI(Q,\alpha)_{-\langle \cdot, \beta \rangle}.$$
\end{lemma}

Now, we can define $(\alpha \circ \beta)_{Q}$ by
$$(\alpha \circ \beta)_{Q}=\dim \SI(Q,\beta)_{\langle \alpha,\cdot
\rangle}=\dim \SI(Q,\alpha)_{-\langle \cdot, \beta \rangle}.$$

(When no confusion arises, we drop the subscript $Q$.)

\subsection{Exceptional sequences}\label{embed} A dimension vector $\beta \in \NN^{Q_0}$ is said to be a \emph{Schur root}
if there exists a $\beta$-dimensional representation $W \in
\Rep(Q,\beta)$ such that $\End_{Q}(W) \cong \CC.$

\begin{definition} A sequence of dimension vectors $\varepsilon_1, \dots, \varepsilon_r$ is called an
exceptional sequence if:
\begin{enumerate}
\renewcommand{\theenumi}{\arabic{enumi}}

\item each $\varepsilon_i$ is a real Schur root, i.e.,
$\varepsilon_i$ is a Schur root and $\langle
\varepsilon_i,\varepsilon_i \rangle =1,$ for all $1 \leq i \leq
r;$

\item $(\varepsilon_i \circ \varepsilon_j)_{Q} \neq 0,$ for all $1
\leq i<j \leq r.$

\end{enumerate}
\end{definition}

The following theorem will be quite useful for us (for a more
general version, see \cite{DW2}):

\begin{theorem} \label{embed-thm}\cite[Theorem 2.39]{DW2} Let $\varepsilon_1,\varepsilon_2$ be an exceptional sequence for a quiver
$Q$ without oriented cycles. Assume that $\langle
\varepsilon_2,\varepsilon_1 \rangle = -l,$ where $l$ is some
non-negative integer. Define a new quiver $\theta(l)$ with set of
vertices $\theta(l)_0=\{1,2\}$ and $l$ arrows from vertex $2$ to
vertex $1.$ Consider the linear transformation
$$
I:\NN^{\theta(l)_0}=\NN^2 \rTo \NN^{Q_0}
$$
defined by
$$
I(\beta_1,\beta_2)=\beta_1\varepsilon_1+\beta_2\varepsilon_2,
$$
for all dimension vectors $\beta=(\beta_1,\beta_2) \in
\NN^{\theta(l)_0}.$

If $\alpha, \beta \in \NN^{\theta(l)_0}$ are so that $(\alpha
\circ \beta)_{\theta(l)} \neq 0$ then
$$
(\alpha \circ \beta)_{\theta(l)}=(I(\alpha) \circ I(\beta))_{Q}.
$$
\end{theorem}

The quiver $\theta(l)$ that appears in Theorem \ref{embed-thm} is
called \emph{the generalized Kronecker quiver.} As we will see in
Section \ref{cex-sec}, this particular quiver will be our main
source of Littlewood-Richardson coefficients. It has been proved
in \cite{DW2} that the map $I$ in the theorem above allows one to
"embed" much of the combinatorics of the "new quiver" $\theta(l)$
into the combinatorics of the original quiver $Q$. For this
reason, we refer to Theorem \ref{embed-thm} as the "embedding
theorem".

\subsection{Polynomiality for semi-invariants and
(non-)log-concavity.}\label{poly-sec} We are interested in how the
dimensions $N\alpha \circ \beta=\dim_{\CC}\SI(Q,\beta)_{N\langle
\alpha, \cdot \rangle}$ and $\alpha \circ
N\beta=\dim_{\CC}\SI(Q,\alpha)_{-N\langle \cdot, \beta \rangle}$
vary as $N \in \ZZ_{\geq 0}$ varies.

\begin{proposition}\cite[Corollary 1]{DW3} Let $\alpha, \beta$ be two dimension
vectors such that $\alpha \circ \beta \neq 0.$ There exist
polynomials $P, Q \in \QQ[X]$ (both depending on $\alpha$ and
$\beta$) with $P(0)=Q(0)=1,$ and
$$
N\alpha \circ \beta=P(N),~\forall{N \geq 0},
$$
and
$$
\alpha \circ N\beta=Q(N),~\forall{N \geq 0}.
$$
\end{proposition}

\begin{remark}\label{asympt-rmk} Note that there is a sufficiently large
integer $N_0>0$ such that $p(t)={P(t+1) \over P(t)}$ and
$q(t)={Q(t+1) \over Q(t)}$ are weakly decreasing functions on
$[N_0,\infty).$ In other words, we have
\begin{eqnarray}\label{ineq1}
((N+1)\alpha \circ \beta) \cdot ((N-1)\alpha \circ \beta) \leq
(N\alpha \circ \beta)^2,
\end{eqnarray}
\begin{eqnarray}
(\alpha \circ (N+1)\beta) \cdot (\alpha \circ (N-1)\beta) \leq
(\alpha \circ N\beta)^2,
\end{eqnarray}
for every $N>N_0.$

Thus the dimensions of spaces of semi-invariants are
\emph{asymptotically} log-concave (in each argument).
\end{remark}

\begin{example}\label{asympt-ex} For an integer $r \geq 1,$ let $T_{r,r,r}$ be the following triple flag quiver
with arms of length $r:$
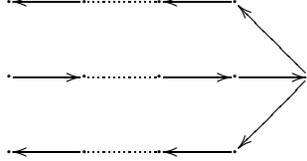
$$
\xy (0, 0)*{\cdot}="a";
        (-10, 10)*{\cdot}="a1";
        (-20,10)*{\cdot}="a2";
        (-30,10)*{\cdot}="a3";
        (-40,10)*{\cdot}="a4";
        (-10,0)*{\cdot}="b1";
        (-20,0)*{\cdot}="b2";
        (-30,0)*{\cdot}="b3";
        (-40,0)*{\cdot}="b4";
        (-10,-10)*{\cdot}="c1";
        (-20,-10)*{\cdot}="c2";
        (-30,-10)*{\cdot}="c3";
        (-40,-10)*{\cdot}="c4";
        {\ar@{->} "a";"a1"};
        {\ar@{->} "a1";"a2"};
        {\ar@{.} "a2";"a3"};
        {\ar@{->} "a3";"a4"};
        {\ar@{->} "b1";"a"};
        {\ar@{->} "b2";"b1"};
        {\ar@{.} "b2";"b3"};
        {\ar@{->} "b4";"b3"};
        {\ar@{->} "a";"c1"};
        {\ar@{->} "c1";"c2"};
        {\ar@{.} "c2";"c3"};
        {\ar@{->} "c3";"c4"};
    \endxy
$$

Now, given a triple $(\lambda, \mu, \nu)$ of partitions of length
at most $r,$ one can construct dimension vectors $\alpha$ and
$\beta$ (see for example \cite{DW1}) such that
$$
N \alpha \circ \beta = c_{N\mu,N\nu}^{N\lambda},
$$
for all $N \geq 1.$ This calculation together with inequality
$(\ref{ineq1})$ shows that the Littlewood-Richardson coefficients
are \emph{asymptotically} log-concave (compare with \cite[Section
3.5]{Oko}).
\end{example}

Next, we are going to show that the log-concavity property for
semi-invariants fails in many cases. For $\beta \in \NN^{Q_0}$ a
dimension vector and $\sigma \in \ZZ^{Q_0}$ a weight of $Q,$ we
define
$$
\sigma(\beta)=\sum_{x \in Q_0}\sigma(x)\beta(x).
$$

\begin{definition}\cite[Proposition 3.1]{K} Let $\beta$ be a dimension vector and
$\sigma$ be a weight such that $\sigma(\beta)=0.$ A
$\beta$-dimensional representation $W \in \Rep(Q,\beta)$ is said
to be:
\begin{enumerate}
\renewcommand{\theenumi}{\arabic{enumi}}

\item \emph{$\sigma$-semi-stable} if $\sigma(\underline d_{W'})
\leq 0$ for every subrepresentation $W'$ of $W;$

\item \emph{$\sigma$-stable} if $\sigma(\underline d_{W'})<0$ for
every proper subrepresentation $0 \neq W' \varsubsetneq W.$
\end{enumerate}
\end{definition}

We say that a dimension vector $\beta$ is
\emph{$\sigma$(-semi)-stable} if there exists a
$\sigma$(-semi)-stable representation $W \in \Rep(Q,\beta).$

Let $\beta$ be a $\sigma$-semi-stable dimension vector. The set of
$\sigma$-semi-stable representations in $\Rep(Q,\beta)$ is denoted
by $\Rep(Q,\beta)^{s.s.}_{\sigma}$ while the set of
$\sigma$-stable representations in $\Rep(Q,\beta)$ is denoted by
$\Rep(Q,\beta)^{s}_{\sigma}.$ The one dimensional torus
$$
T=\{(t\Id_{\beta(x)})_{x \in Q_0} \mid t \in \CC^{*} \} \subseteq
\GL(\beta)
$$
acts trivially on $\Rep(Q,\beta)$ and so there is a well-defined
action of $\PGL(\beta)={\GL(\beta)/T}$ on $\Rep(Q,\beta).$ Using
methods from geometric invariant theory, one can construct the
following GIT-quotient of $\Rep(Q,\beta)$:
$$
\mathcal M (Q,\beta)^{s.s.}_{\sigma}=\Proj (\oplus_{n \geq
0}\SI(Q, \beta)_{n\sigma}).
$$
It was proved by King \cite{K} that $\mathcal M
(Q,\beta)^{s.s.}_{\sigma}$ is a categorical quotient of
$\Rep(Q,\beta)^{s.s.}_{\sigma}$ by $\PGL(\beta).$ Note that
$\mathcal M(Q,\beta)^{s.s.}_{\sigma}$ is an irreducible projective
variety, called the moduli space of $\beta$-dimensional
$\sigma$-semi-stable representations (for more details, see
\cite{K}).

For the remainder of this section, we assume that $\beta$ is a
$\sigma$-stable dimension vector. Then there is a non-empty open
subset $\mathcal M(Q,\beta)^s_{\sigma} \subseteq \mathcal
M(Q,\beta)^{s.s.}_{\sigma}$ which is a geometric quotient of
$\Rep(Q,\beta)^s_{\sigma}$ by $\PGL(\beta).$ Now, a
$\sigma$-stable representation must be a Schur representation and
so its stabilizer in $\PGL(\beta)$ is zero dimensional. It follows
that
$$
\dim \mathcal M (Q,\beta)^{s.s.}_{\sigma}=1-\langle \beta,\beta
\rangle.
$$

Let us further assume that $\langle \beta, \beta \rangle <0$ (that
is to say, $\beta$ is \emph{imaginary and non-isotropic}). Then it
is known that $m\beta$ stays $\sigma$-stable (see for example
\cite[Proposition 3.16]{DW2}) and hence
$$
\dim \mathcal M(Q,m\beta)^{s.s.}_{\sigma}=1-m^2 \langle \beta,
\beta \rangle,
$$
for every integer $m \geq 1.$ Now, write $\sigma = \langle \alpha,
\cdot \rangle$ for some dimension vector $\alpha.$ If we fix $m$
then $n\alpha \circ m\beta$ has degree $1-m^2\langle \beta, \beta
\rangle$ as a polynomial in $n.$ Therefore, when $n>0$ is
sufficiently large, we must have
\begin{eqnarray}\label{cex-semi-inv}
n\alpha \circ 2\beta > (n\alpha \circ \beta)^2.
\end{eqnarray}
Indeed, the left hand-side of the above inequality is a polynomial
in $n$ of degree $1-4\langle \beta, \beta \rangle$ while the right
hand-side is a polynomial of degree $2-2\langle \beta, \beta
\rangle$ and $1-4\langle \beta, \beta \rangle > 2-2\langle \beta,
\beta \rangle.$

Note that inequality $(\ref{cex-semi-inv})$ gives counterexamples
to the log-concavity property for semi-invariants.

\section{Counterexamples}\label{cex-sec}
In this section, we first give a different proof of Theorem
\ref{prop-cex} and then present more counterexamples. In
particular, we provide counterexamples to Kirillov's $q$-Log
concavity conjecture for parabolic Kostka polynomials (see
Proposition \ref{parK-cex-prop}).

\subsection{Littlewood-Richardson coefficients from star and
generalized Kronecker quivers}\label{compute}

It is well-known that the Littlewood-Richardson coefficients can
be viewed as dimensions of spaces of semi-invariants of star
quivers (see for example \cite{CC1}, \cite{DW1}). Now, let us
consider the star quiver $T_{4,3,4}$ with the following
orientation:

$$
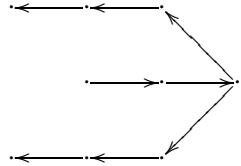
\xy (0, 0)*{\cdot}="a";
        (-10, 10)*{\cdot}="b";
        (-20,10)*{\cdot}="c";
        (-30,10)*{\cdot}="d";
        (-20,0)*{\cdot}="e";
        (-10,0)*{\cdot}="f";
        (-30,-10)*{\cdot}="g";
        (-20,-10)*{\cdot}="h";
        (-10,-10)*{\cdot}="i";
        {\ar@{->} "a";"b"};
        {\ar@{->} "b";"c"};
        {\ar@{->} "c";"d"};
        {\ar@{->} "f";"a"};
        {\ar@{->} "e";"f"};
        {\ar@{->} "a";"i"};
        {\ar@{->} "i";"h"};
        {\ar@{->} "h";"g"};
    \endxy
$$

We are going to reduce the problem of computing semi-invariants of
$T_{4,3,4}$ to that of computing semi-invariants of (rather small)
generalized Kronecker quivers.

Let us recall that for every integer $n \geq 1,$ we define

$$
\lambda(n)=(4^n,3^{2n},2^n),\text{~and~} \mu(n)=(3^n,2^n,1^n).
$$

\begin{proposition} \label{embed-prop}Let $\theta(3)$ be the generalized Kronecker quiver with $3$
arrows and vertices labelled $1,2:$
$$
\theta(3):~ \xy     (0,0)*{1}="a";
        (10,0)*{2}="b";
        {\ar@3{<-} "a";"b" };
    \endxy
$$

Then
$$
\dim \SI(\theta(3),(n,n))_{(-m,m)}=c^{m\lambda(n)}_{m \mu(n),m
\mu(n)},
$$
for every $m, n \geq 1.$
\end{proposition}

\begin{proof} Let us consider the exceptional sequence of $T_{4,3,4}$
given by:
$$
\varepsilon_1=
\begin{matrix}
1&2&3& \\
&0&3&4,\\
1&2&3&
\end{matrix}
$$
and
$$
\varepsilon_2=
\begin{matrix}
0&0&0& \\
&1&0&0.\\
0&0&0&
\end{matrix}
$$

Since $\langle \varepsilon_2, \varepsilon_1 \rangle =-3,$ we know
that the generalized Kronecker quiver $\theta(3)$ can be embedded
in $T_{4,3,4}$ by Theorem \ref{embed-thm}. In particular, if
$\alpha=(n,n)$ and $\beta=(m,2m)$ are dimension vectors for
$\theta(3)$ then
$$
(\alpha \circ \beta)_{\theta(3)} =(I(\alpha) \circ
I(\beta))_{T_{4,3,4}},
$$
where
$$
I(\alpha)=
\begin{matrix}
n&2n&3n& \\
&n&3n&4n,\\
n&2n&3n&
\end{matrix}
$$
and
$$
I(\beta)=
\begin{matrix}
m&2m&3m& \\
&2m&3m&4m.\\
m&2m&3m&
\end{matrix}
$$

Next, computing with Schur functors (see \cite{CC1} or \cite{DW1}
for explicit computations) we obtain
$$
I(\alpha) \circ I(\beta)=c^{m\lambda(n)}_{m \mu(n),m \mu(n)},
$$
and so we are done.
\end{proof}

\begin{proof}[Another proof of Theorem \ref{prop-cex}]
By Proposition \ref{embed-prop}, we only need to compute the
dimensions of the spaces $\SI(\theta(3),(n,n))_{(-m,m)}$ when
$m=1,2.$ For this, we first decompose the affine coordinate ring
of $\Rep(\theta(3), (n,n))$ as a  direct sum in which the summands
are tensor products of irreducible representations of $GL(n)$'s.
For convenience, let us write $V=\CC^n, W=\CC^n.$ Then we have
$$
\begin{aligned}
\CC[\Rep(\theta(3), (n,n))]&=\CC[\Hom(W,V)\oplus \Hom(W,V) \oplus
\Hom(W,V)]\\
&=S(W \otimes V^*)\otimes S(W \otimes V^*)\otimes S(W \otimes
V^*).
\end{aligned}
$$
Using Cauchy's formula \cite[page 121]{F3}, we obtain that
$$
S(W \otimes V^*)=\oplus S^{\mu}W \otimes S^{\mu}V^*
$$
as $\GL(V) \times \GL(W)$-modules, where the sum is over all
partitions $\mu$ with at most $n$ non-zero parts. Hence, we have
$$
\begin{aligned}
&\CC[\Rep(\theta(3), (n,n))]^{\SL(V) \times \SL(W)}=\\
&=\oplus (S^{\mu(1)}V^* \otimes S^{\mu(2)}V^* \otimes
S^{\mu(3)}V^*)^{\SL(V)} \otimes (S^{\mu(1)}W \otimes S^{\mu(2)}W
\otimes S^{\mu(3)}W)^{\SL(W)},
\end{aligned}
$$
where the sum is over all partitions $\mu(1), \mu(2), \mu(3)$ with
at most $n$ non-zero parts. Sorting out those semi-invariants of
weight $(-m,m),$ it is easy to see that
$$\SI(\theta(3),(n,n))_{(-m,m)}= \oplus
({\det}_{V}^{m}\otimes \otimes_{i=1}^3
S^{\mu(i)}V^*)^{\GL(V)}\otimes ({\det}_{W}^{-m} \otimes
\otimes_{i=1}^3 S^{\mu(i)}W )^{\GL(W)},
$$
where the sum is over all partitions $\mu(1), \mu(2), \mu(3)$ with
at most $n$ non-zero parts. For our purposes it is useful to work
with conjugate partitions in the identity above. So, we can write
\begin{eqnarray}\label{eqn0}
\dim \SI(\theta(3),(n,n))_{(-m,m)}=\sum \left( c_{\lambda(1),
\lambda(2), \lambda(3)}^{(n^m)} \right)^2,
\end{eqnarray}
where the sum is over all partitions $\lambda(1), \lambda(2),
\lambda(3)$ (with at most $m$ non-zero parts).

Next, it is easy to see that
\begin{eqnarray}\label{eqn1}
c_{\lambda(1), \lambda(2), \lambda(3)}^{(n^m)} \leq 1,
\end{eqnarray}
for $m \in \{1,2 \}.$ Indeed, one can either check this directly
with the Littlewood-Richardson rule or view these coefficients as
dimensions of spaces of semi-invariants for a quiver of type
$\mathbb D_4.$

Therefore, $\dim \SI(\theta(3),(n,n))_{(-1,1)}$ is simply the
number of monomials in three (commuting) variables of degree $n,$
and so
$$
c_{\mu(n),\mu(n)}^{\lambda(n)}=\binom{n+2}{2}.
$$

Now, let $\lambda(i)=(\lambda_1(i),\lambda_2(i)), ~1 \leq i \leq
3$ be three partitions with at most two non-zero parts. We claim
that
\begin{gather}
|\lambda(1)|+|\lambda(2)|+|\lambda(3)|=2n, \label{eqn2}\\
n-\lambda_1(i)-\lambda_2(j)-\lambda_2(k) \geq 0, \text{~where~}
\{i,j,k\}=\{1,2,3 \},\label{eqn3}
\end{gather}
give a (minimal) list of necessary and sufficient Horn
inequalities for the non-vanishing of the Littlewood-Richardson
coefficient $c_{\lambda(1), \lambda(2), \lambda(3)}^{(n^2)}.$ This
follows from \cite[Theorem 17]{F1}. Alternatively, one can deduce
this claim from the description of the so called cone of effective
weights for a type $\mathbb D_4$ quiver.

From $(\ref{eqn0})-(\ref{eqn3}),$ we obtain that $\dim
\SI(\theta(3),(n,n))_{(-2,2)}$ equals the cardinality of the set
$\mathcal S$ of all triples $(\lambda(1),\lambda(2), \lambda(3))$
of partitions with at most $2$ non-zero parts satisfying the
conditions $(\ref{eqn2})$ and $(\ref{eqn3}).$

Note that every $(\lambda(1), \lambda(2), \lambda(3)) \in \mathcal
S$ gives rise to a monomial $X_1^{n_1} \cdot X_2^{n_2}\cdot
X_3^{n_3} \cdot X_4^{n_4} \cdot X_5^{n_5} \cdot X_6^{n_6}$ of
degree $n,$ where
\begin{align*}
n_1&=n-\lambda_1(1)-\lambda_2(2)-\lambda_2(3)&n_2&=\lambda_2(1)\\
n_3&=n-\lambda_1(2)-\lambda_2(3)-\lambda_2(1)&n_4&=\lambda_2(2)\\
n_5&=n-\lambda_1(3)-\lambda_2(1)-\lambda_2(2)&n_6&=\lambda_2(3).
\end{align*}

It is clear that in this way we get a bijection from $\mathcal S$
to the set of all monomials in six (commuting) variables of degree
$n.$ So, we have
$$
\dim \SI(\theta(3),(n,n))_{(-2,2)}=\binom{n+5}{5}
$$
and this finishes the proof.
\end{proof}

\begin{remark} It is worth pointing out that using the same ideas
as above one can construct non-log-concave Littlewood-Richardson
coefficients for every star quiver $T_{p,q,r}$ of wild
representation type.
\end{remark}

\subsection{Non-log-concave parabolic Kostka numbers} \label{parKostka} In this section, we consider some rather special Littlewood-Richardson coefficients.
Let $\lambda$ be a partition and let $R=((m_1^{l_1}), \dots,
(m_k^{l_k}))$ be a sequence of rectangular partitions. Then the
parabolic Kostka number $K_{\lambda,R}$ associated to $\lambda$
and $R$ is defined by
$$
K_{\lambda,R}=\dim_{\CC}(S^{\lambda}(V)^* \otimes
S^{(m_1^{l_1})}(V)\otimes \cdots \otimes
S^{(m_k^{l_k})}(V))^{\GL(V)},
$$
where $V$ is a complex vector space of sufficiently large
dimension. In general, it is well-known that $K_{\lambda,R}$ is
the value at $q=1$ of the corresponding parabolic Kostka
polynomial (see \cite[Chapter 4]{Kir} and the reference therein).

If $R=((m_1^{l_1}), \dots, (m_k^{l_k}))$ is a sequence of
rectangles and $N \geq 1$ is an integer, we define $NR$ to be the
sequence of rectangles $NR=(((Nm_1)^{l_1}),\dots,
((Nm_k)^{l_k})).$ The log-concavity conjecture for parabolic
Kostka numbers (compare with the more general version
\cite[Conjecture 6.17]{Kir}) is:

\begin{conjecture} \label{parK-conj} Let $\lambda$ be a partition and $R$ be a
sequence of rectangular partitions. Then
$$
K_{(N+1)\lambda,(N+1)R} \cdot K_{(N-1)\lambda,(N-1)R} \leq
(K_{N\lambda,N R})^2,
$$
for every integer $N \geq 1.$
\end{conjecture}

Our next proposition shows that Conjecture \ref{parK-conj} fails
in general.

\begin{proposition} \label{parK-cex-prop} For every $n \geq 1,$
consider
$$
\lambda(n)=(2^n,1^{2n})
$$
and
$$
R(n)=((1^n),(1^n),(1^n),(1^n)).
$$
Then
$$
K_{\lambda(n), R(n)}=\binom{n+2}{2} \text{~and~} K_{2\lambda(n),
2R(n)}= \binom{n+5}{5}.
$$

Consequently, when $n \geq 21,$ Conjecture \ref{parK-conj} fails
for $\lambda=\lambda(n),$ $R=R(n)$ and $N=1.$
\end{proposition}

\begin{proof}
To obtain parabolic Kostka numbers, we work with the following
star quiver $Q:$
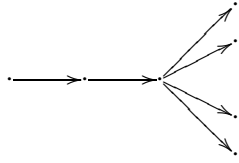
$$
\xy (-20, 0)*{\cdot}="a";
        (-10, 0)*{\cdot}="b";
        (0,0)*{\cdot}="c";
        (10,10)*{\cdot}="d";
        (10,5)*{\cdot}="e";
        (10,-5)*{\cdot}="f";
        (10,-10)*{\cdot}="g";
        {\ar "a";"b"};
        {\ar "b";"c"};
        {\ar "c";"d"};
        {\ar "c";"e"};
        {\ar "c";"f"};
        {\ar "c";"g"};
    \endxy
$$

Let $\varepsilon_1, \varepsilon_2$ be the following exceptional
sequence of $Q:$

$$
\varepsilon_1=
\begin{matrix}
&  & & &1\\
 &  & & &1 \\
&0&3&4& &,\\
&  & & &1\\
&  & & &1
\end{matrix}
$$
and
$$
\varepsilon_2=
\begin{matrix}
&  & & &0\\
 &  & & &0 \\
&1&0&0& &.\\
&  & & &0\\
&  & & &0
\end{matrix}
$$

Reasoning as in Proposition \ref{embed-prop}, we get
$$
\dim \SI(\theta(3),(n,n))_{(-m,m)}=K_{m\lambda(n),m R(n)},
$$
for every $m, n \geq 1.$ The proof follows from that of Theorem
\ref{prop-cex}.
\end{proof}

\begin{remark} Note that the parabolic Kostka
numbers appearing in Proposition \ref{parK-cex-prop} can be
written as Littlewood-Richardson coefficients:
$$
K_{m\lambda(n),
mR(n)}=c_{((3m)^n,(2m)^n,m^n),((2m)^n,m^{2n})}^{((4m)^n,(3m)^n,(2m)^n,m^n)},
$$
for every integer $m \geq 1.$ Indeed, this follows immediately
from \cite[Proposition 9]{Z}.
\end{remark}

\section*{Acknowledgment} We would like to thank the anonymous referee for useful comments that
helped improve the exposition of the paper .

%\bibliography{biblio}

\end{document}